\documentclass{ifacconf}
\usepackage{graphicx}      
\usepackage{natbib}        
\newtheorem{theorem}{Theorem}

\newtheorem{remark}{Remark}

\newtheorem{assumption}{Assumption}
\usepackage{amssymb}
\usepackage{amsmath}
\usepackage{mathrsfs}
\usepackage{makecell}
\usepackage{algorithm}
\usepackage{algorithmic}
\usepackage{xcolor}
\usepackage{color}
\usepackage{pgfplots}
\usepackage{theorem}
\usepackage{enumitem}
\usepackage{empheq}
\usepackage{ulem}
\usepackage[thicklines]{cancel}
\makeatletter
\let\old@ssect\@ssect 
\makeatother
\usepackage[hidelinks]{hyperref}
\usepackage{soul}
\makeatletter
\def\@ssect#1#2#3#4#5#6{%
	\NR@gettitle{#6}
	\old@ssect{#1}{#2}{#3}{#4}{#5}{#6}
}
\makeatother

\setitemize[1]{itemsep=5pt,partopsep=3pt,parsep=\parskip,topsep=5pt}

\newcommand{\imag}{\mathfrak{i}}
\newcommand{\tr}{\top}

\newcommand{\Shijie}[1]{#1\color{black}}

\newcommand{\newshijie}[1]{\color{black}#1\color{black}}

\newcommand{\Proof}{\noindent \textbf{Proof.}$\;\;$}
\newcommand{\bracketcite}[1]{ (\cite{1})}

\def\equationautorefname~#1\null{Equation~(#1)\null}

\begin{document}
	\begin{frontmatter}
		
		\title{Optimal 
Demand Shut-offs
			of AC Microgrid\\ using AO-SBQP Method} 

		
		\thanks{All five authors contribute equally. 
		}
		
		\author[First,Second]{Xu Du} 
		\author[First,Third]{Shijie Zhu} 
		\author[First,Third]{Yifei Wang}
		\author[First,Third]{Boyu Han}  
		\author[First,Third]{Xiaohe He}

		\address[First]{The School of Information Science and Technology, ShanghaiTech University, Shanghai, China}
		\address[Second]{National Key Laboratory  on Wireless Communications, University of Electronic Science and Technology of China, \hspace{-0.1cm}Chengdu,\hspace{-0.1cm} China\\ (e-mail: duxu@uestc.edu.cn).}
		\address[Third]{
			Innovation Academy for Microsatellites, Chinese Academy of Sciences, Shanghai, China\\
		University of Chinese Academy of Sciences, Beijing, China\\ (e-mail: zhushj, wangyf6, hanby, hexh@shanghaitech.edu.cn)}
	
		
		\begin{abstract}
			Microgrids are increasingly being utilized to improve the resilience and operational flexibility of power grids, and act as a backup power source during grid outages. 
			However, it necessitates that the microgrid itself could provide power to the critical loads. This paper presents an algorithm named alternating optimization based sequential boolean quadratic programming tailored for solving optimal demand shut-offs problems arising in microgrids. 
			Moreover, we establish local superlinear convergence of the proposed approximate Boolean quadratic programming method over nonconvex problems. In the end, the performance of the proposed method is illustrated on the modified IEEE 30-bus case study.

		\end{abstract}
		
		
	\end{frontmatter}
	
	\section{Introduction}\label{sec: introduction}
	
	Microgrids are small-scale power grids that can operate autonomously or in collaboration with other small power grids, which are increasingly being utilized to increase the resilience and operational flexibility of power grids. They act as a backup power supply in the case of grid outages caused by devastating disasters.
	\cite{abbey2014powering} summarized the measures taken by Sendai region in Japan to cope with the shortage of circuit supply caused by the nuclear power accident after the 2011 East Japan Earthquake. Similarly,  \cite{panora2014real} describes a successful case of rapid restoration of local microgrid integrated energy systems after infrastructure damage caused by Superstorm Sandy in Manhattan Island, 2012. 
	However, this necessitates that each isolated microgrid itself be resistant and formulate its own optimal power supply strategy with the shortage of energy supply, which is still a challenging problem.
	 Basically, a potential solution is to abstract the above power grid optimization problem in a mixed Boolean nonlinear programming fashion (MBNLP).

	To the best of our knowledge, classical optimization in power grid consists optimal power flow (OPF,\;\cite{frank2012optimal1}), optimal reactive power dispatch (\cite{zhu2015optimization}), power system state and parameter estimation (\cite{Monticelli1999,du2019distributed}). Recently, \cite{du2020optimal,du2022approximations} proposes the method of optimal experimental design (OED) in order to extract more relative information for assisting the admittance estimation process. Moreover, \cite{Du2021} offers an adaptive method for balancing OED and the OPF cost.
	 These mentioned power grid optimization problems are smooth and can be solved directly with interior point method, while notable recent researches optimize the discrete decision variables at the same time. 
	 \cite{rhodes2020balancing} and \cite{kody2022sharing} modeled the direct current (DC) optimal power shut-in problems in a mixed integer linear program (MILP) framework and solved them with Gurobi. However, only few literature can be found for nonconvex MBNLP in the area of alternating current (AC) power systems. On the other hand, from the algorithmic level, the solver is based on branch and bound (B$\&$B) method (\cite{morrison2016branch}), which needs to establish a tree storage structure to explore each integer variable with low efficiency. \cite{luo2010sdp} solves Boolean optimization in a semi-definite relaxation (SDR) fashion, however, with matrix variables. 
	  Solving nonconvex MBNLP accurately and efficiently remains an open problem in our view.
	
	Recently, a quadratic programming with linear complementary constraint (LCQP) problem is well studied by a series of literatures (\cite{lcqp,hall2022lcqpow}). 
	By sequentially solving the QP problem with the corresponding linearlized penalty term, the complementary constraint can be reached with finite steps. Inspired by the above literatures, alternating optimization based sequential boolean quadratic
	programming (AO-SBQP) (\cite{2022Alternating}) is proposed to set up a bridge between LCQP and MBNLP, leading the solution of optimization problems with Boolean variables no longer rely on B$\&$B method with tree structure searching or SDR with matrix variables.
	
	In this paper, the idea of AO-SBQP method is inherited, and the algorithm is modified in the specific optimal distribution of limited supply scenario in microgrid. Unlike \cite{2022Alternating}, our considered model is still a non-convex problem except the Boolean variable constraints. Moreover, a local convergence analysis for the approximate BQP with constraints is introduced.

	The rest of this paper is organized as follows: Section
	\ref{sec: power system} reviews the basics concepts of power grid.
	Sections \ref{sec:model} proposes the optimal distribution of limited supply model. Section \ref{sec:LCQP} describes the AO-SBQP algorithm in detail. And the numerical result is shown in Section \ref{sec: numerical result}.

	\textit{Notation:}  For $a \in \mathbb{R}^n$ and $\mathcal{C}\subseteq\{1,...,n\}$, ${[a_i]_{i \in \mathcal C}\in\mathbb{R}^{|\mathcal{C}|}}$ collects all components of $a$ whose index $i$ is in $\mathcal C$.  Similarly, for  $A \in \mathbb R^{n \times l}$ and $\mathcal S \subseteq \{ 1, \ldots, n \} \times \{ 1, \ldots, l \}$, $[A_{i,j}]_{(i,j) \in \mathcal S}\in\mathbb{R}^{|\mathcal{S}|}$ denotes the {concatenation of} $A_{i,j}$ for all $(i,j)\in\mathcal{S}$. $\imag = \sqrt{-1}$ denotes the imaginary unit, such that  $\mathrm{Re}(z) +\imag\cdot \mathrm{Im}(z)=z \in \mathbb C$, and $\hat a$ denotes the estimated value of $a$. Moreover, $\mathbf{1}$ represents a vector with all entries being one. 
	
	\section{AC Power Grid Model}\label{sec: power system}
\begin{figure}[H] 
		\centering
		\includegraphics[width=0.6\linewidth]{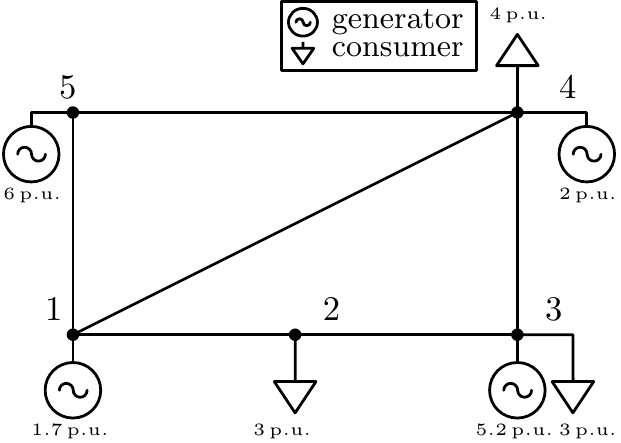}
		\caption{Modiﬁed IEEE 5-bus system from \cite{Li2010} with 4 generators and 3 consumers.} 
		\label{fig:ieee5bus}
	\end{figure}

	Consider a power grid defined by the triple $(\mathcal{N},\mathcal{L},(G+\imag B))$, where $\mathcal{N} = \{1,2,\dots N\}$ represents the set of buses, $\mathcal{L}\subseteq \mathcal{N}\times \mathcal{N}$ denotes  transmission lines and $(G+\imag B)\in \mathbb{C}^{N\times N}$ is the complex and potentially sparse Laplacian admittance matrix
	\[
	(G_{k,l}+\imag B_{k,l})\doteq\left\{
	\begin{array}{ll}
		\sum\limits_{i \neq k} \left( g_{k,i} + \imag \, b_{k,i} \right) & \text{if} \; k=l, \\  [0.25cm]
		- \left( g_{k,l} + \imag \, b_{k,l} \right) & \text{if} \; k \neq l.
	\end{array}
	\right.
	\]
	Here, $g_{k,l}$ and $b_{k,l}$ are the conductances and susceptances of the transmission line $(k,l)\in \mathcal{L}$, which aims to connect the buses.
	Note that $(G_{k,l}+\imag B_{k,l})=0$ if $(k,l)\notin \mathcal{L}$.
	The set  $\mathcal G\subseteq \mathcal N$ collects all nodes equipped with generators and $\mathcal D\subseteq \mathcal N$ collects all nodes with power demands.
	Figure \ref{fig:ieee5bus} shows a 5-bus system with $\mathcal{N}=\{1,\dots,5\}$, $\mathcal G = \{1,3,4,5\}$, $\mathcal D = \{2,3,4\}$.

	Let $v_k$ denote the voltage amplitude at the $k$-th node and $\theta_{k}$ the corresponding voltage angle, $\theta_{k,l}$ denotes the angle difference between node $k$ and $l$. Throughout this paper, we assume that the voltage magnitude and the voltage angle at the first node (the slack node) are fixed, $\theta_1 = 0 \quad \text{and} \quad v_1 = \mathrm{const} \; 
	$.  We refer
	to \cite{Du2021} and \cite{du2022approximations} for further discussion.
	The state variables of the system is defined as
	\begin{equation*}
		x\doteq 
		\begin{bmatrix}
			v_k,
			\theta_k
		\end{bmatrix}_{k \in \mathcal N}^\top \in\mathbb{R}^{2|\mathcal N|}.\; 
	\end{equation*}
	Moreover, we have active and reactive power generation of generators $p_k^g$ and $q_k^g$ for all $k\in \mathcal G$, and  $d_k = (p_k^d ,\;q_k^d)^\top$ denote the active and reactive power demand at demand nodes $k\in \mathcal D $.
	We consider active and reactive power supply of all  generators
	\begin{equation*}
		u\doteq 
		\begin{bmatrix}
			p_k^g,
			q_k^g
		\end{bmatrix}_{k \in \mathcal G }^\top\in\mathbb{R}^{2|\mathcal G|} \; 
	\end{equation*}
	as the system input variables. 
	
	The active and reactive power flow over the transmission line $(k,l) \in \mathcal L$ is given by
	\begin{align*}
		\hspace{-0.5cm}\Pi_{k,l}(x) \doteq &v_k^2 \begin{bmatrix}
			g_{k,l} \\[0.16cm]
			- b_{k,l}
		\end{bmatrix}
		\\& - v_kv_l\begin{bmatrix}
			g_{k,l},  & b_{k,l} \\[0.16cm]
			-b_{k,l}, & g_{k,l} 
		\end{bmatrix}
		\hspace{-0.05cm}
		\hspace{-0.05cm}
		\hspace{-0.05cm}
		\begin{bmatrix}
\cos(\theta_{k,l}) \\[0.16cm]
			\sin(\theta_{k,l})
		\end{bmatrix}
		\hspace{-0.05cm}.
	\end{align*}

	The total power outflow  from node $k \in \mathcal{N}$ is given by
\begin{equation}
\begin{split}
P_k(x) \doteq& v_k^2 \sum_{l \in \mathcal N_k}  
\begin{bmatrix}
g_{k,l} \\[0.16cm]
-b_{k,l}
\end{bmatrix} \notag \\
& - v_k \sum_{l \in \mathcal N_k} v_l
\begin{bmatrix}
g_{k,l},  & b_{k,l} \\[0.16cm]
-b_{k,l}, & g_{k,l}
\end{bmatrix}
\begin{bmatrix}
\cos(\theta_{k,l}) \\[0.16cm]
\sin(\theta_{k,l})
\end{bmatrix} \\
\doteq&\sum_{l\in \mathcal N_k} \Pi_{k,l}(x),
\end{split}
\end{equation}
	where $$\mathcal N_k \doteq \{ l \in \mathcal N \, \mid \, (k,l) \in \mathcal L \, \}$$
	denotes the set of neighbors of node $k \in \mathcal{N}$.
	
	\section{Optimal Distribution of Limited Supply Power Network}\label{sec:model}
	
	In the limited power supply scenario, potentially, \[\sum_{k\in\mathcal{G}}\overline{p_k^g}<\sum_{k\in\mathcal{D}}p_k^d,\]
	here, $\overline{p_k^g}$ denotes the upper bound of active generator power input of bus $k$.
	This leads to a fact that not all the power demand will enjoy stable energy supply. Thus we introduce $y\in \mathbb R^{|\mathcal D|}$ with $y_k\in \{0,1\}$ as an auxiliary switch variable. Then the power supply at a given bus can be expressed as
	\begin{equation*}
		S_k(u,y)\doteq\left\{\begin{aligned}
			&\begin{bmatrix}
				p_k^g,
				q_k^g
			\end{bmatrix}^\top-y^2_k\begin{bmatrix}
				p_k^d,
				q_k^d
			\end{bmatrix}^\top, &k &\in \mathcal D, \mathcal G\\[0.16cm]
			&\begin{bmatrix}
				p_k^g,
				q_k^g
			\end{bmatrix}^\top , \qquad &k &\notin \mathcal D,k\in \mathcal G 
			\\[0.16cm]
			&-y_k^2\begin{bmatrix}
				p_k^d,
				q_k^d
			\end{bmatrix}^\top , \qquad &k &\in \mathcal D,k\notin \mathcal G.
		\end{aligned}\right.
	\end{equation*}
	Thus, the power flow equations can be written in the form 
	\begin{equation}
		\label{eq::powerEq}
		\begin{split}
			&P(x) = S(u,y) \; 
		\end{split},
	\end{equation}
	where
	\begin{align*}
		P(x) &\doteq \left[ P_1(x)^\tr, \ldots, P_N(x)^\tr \right]^\tr,  \\
		\;\;	 S(u,y)&\doteq \left[ S_1(u,y)^\tr, \ldots, S_N(u,y)^\tr \right]^\tr. \; \notag
	\end{align*}
	Notice that 
	$\mathrm{dim}(P) = \mathrm{dim}(x) = 2|\mathcal N| \;.$
	
	Optimal power demand shut-offs (weighted maximum load delivery, \cite{rhodes2020balancing}) can be formulated in the following way for a given positive rank $r_k\in \mathbb R^+$ of each power demand. 
	\begin{subequations}
		\label{eq:power supply}
		\begin{align} 
			\max_{ x, u, y}\;\;& f( x, u,  y)=\sum_{k\in\mathcal{D}}y_k r_k  p^d_k \label{eq: original obj} \\ 
			\quad\mathrm{s.t.}\;\;& P(x) - S(u,y)=0 \\
			&\underline{x}\leq x \leq \overline{x}\\
			&\underline{u}\leq u \leq \overline{u}\\
			&y_k\in\{0,1\}.
		\end{align}
	\end{subequations}

	Here we are trying to provide stable power supply for the \emph{most important consumers}  with limited resource. Unlike other classical optimization problems over power grid (\cite{zhu2015optimization}, \cite{frank2012optimal1}), the Boolean variable $y$ leads Equation~\eqref{eq:power supply} into a non-smooth non-convex MBNLP form. This leads it difficult for the mainstream algorithms such as \emph{interior point method} to be applied directly.
	However, it can be reformulate as a  mathematical programs with \emph{complementarity constraints problem}(MPCC) (\cite{lcqp}, \cite{Hall2021a}). 
	
	We reformulate Problem~\eqref{eq:power supply} in the following form:
	\begin{subequations}\label{eq:object2}
		\begin{align}
			\max_{ x, u, y}\;& E(x,u,y) \label{eq:object}\\ \quad\mathrm{s.t.}\;\;
					&C(x,u,y)\leq 0\label{eq:all equaltity}\\
			&0\leq y_k \perp(1-y_k)\geq 0,\label{eq: MPCC}
		\end{align}
		\end{subequations}
	with 
\begin{equation*}
\begin{split}
E(x,u,y)\doteq\sum_{k\in\mathcal{D}}&y_kr_kp_k^g-\\ &\hspace{-0.35cm}\left(\hspace{-0.1cm}y_kr_kv_k\hspace{-1mm}\sum_{l\in\mathcal{N}}\hspace{-0.1cm} v_l(G_{k,l}\cos(\theta_{k,l})\hspace{-0.07cm}+\hspace{-0.07cm}B_{k,l}\sin(\theta_{k,l}))\hspace{-0.1cm}\right),
\end{split}
\end{equation*}
		\begin{equation}\label{eq::physical constyraints}
	\begin{split}
	C(x,u,y)\doteq\begin{bmatrix}
	&P(x) - S(u,y) \\
		&S(u,y) -P(x) \\
	&\underline{x}- x\\
	&x-\overline{x}\\
	&\underline{u}- u\\
	&u- \overline{u}
	\end{bmatrix}
	\end{split}
	\end{equation}
and
 \begin{subequations}
	\label{cons:lincompls}
	\begin{empheq}[left ={0\leq y_k\perp (1-y_k)\geq 0\Leftrightarrow \empheqlbrace}]{align}
		0&\leq y_k\label{cons:leftcompcons} \\
		0&\leq 1-y_k\label{cons:rightcompcons} \\
		0&=y_k\cdot (1-y_k).\label{cons:bilinearcompcons}
	\end{empheq}
\end{subequations}
Note that the problem is based on AC static model at a specific time point without dynamic. The renewables and storage for microgrid is considered as the part of power supply in the model so considering these modules separately is not essential.
	Apart from the MPCC constraints, Problem~\eqref{eq:object2} is still a nonlinear programming (NLP) problem which is relatively not easy to handle. In the next section, we will introduce a new method called \emph{Alternating Optimization Based Sequential Boolean Quadratic Programming Method (AO-SBQP)}  (\cite{2022Alternating}) that can deal with it.
	
		\section{Alternating Optimization Based Sequential Boolean Quadratic Programming}\label{sec:LCQP}
		
		This section reviews the basic AO-SBQP structure (\cite{2022Alternating}) which solves Problem~\eqref{eq:object2} in a sequential way by optimizing the continuous and Boolean variables into different steps.  
		
		\subsection{Approximation of BQP}
		
		The Lagrangian function of \autoref{eq:object2} (\cite{bertsekas1997nonlinear}) without \eqref{cons:lincompls} is 
		\begin{equation*}
		\begin{split}
&\mathcal L_0(x,u,y,\lambda\Shijie{_0})\doteq E(x,u,y)-\lambda_0^\top C(x,u,y).		\end{split}
		\end{equation*}
		According to the standard penalty reformulation, 		\[\phi(y)\doteq y^\tr(\mathbf{1}-y)\]
		is considered as the bi-linear complementarity penalty function relates to \autoref{cons:bilinearcompcons}.
		A second-order Taylor expansion of $\mathcal L_0(y,\lambda\Shijie{_0})$ (\cite{lcqp}) on penalty with respect to $y$ is
		\begin{equation}\label{eq:approximate}
			\upsilon(y)\doteq\frac{1}{2}y^\tr Q y +(g-\rho \nabla\phi(\tilde y))^\tr y.
		\end{equation}
			Here $0\succ Q\doteq\nabla^2_y\mathcal{L}_0(\tilde x,\tilde u,\tilde y,\tilde{\lambda}\Shijie{_0})\in\mathbb{R}^{|\mathcal D|\times|\mathcal D|}$, $g= \nabla_yE(\tilde x,\tilde u,\tilde y) \in\mathbb{R}^{|\mathcal D|}$, $\rho>0$,
		and $(\tilde x,\tilde u,\tilde y,\tilde{\lambda})$ represents the value of $(x,u,y,\lambda)$ from the last NLP iteration.

		For each iteration, the following simplified QP (\ref{eq:BQP}) needs to be solved, here
		\begin{subequations}\label{eq:BQP}
			\begin{align}
			\max_{y}\;\;& \upsilon(y-\tilde{y})\label{eq:BQP:obj} \\ \quad\mathrm{s.t.}\;\;& b+A\cdot(y-\tilde{y})\geq0 \hspace{-0.5mm}\quad&|\lambda\label{eq:BQP:affcons}\\
			&y\geq0 \;\quad\qquad&|\mu\label{eq: BL}\\ 
			&\mathbf{1}-y\geq0\;\;\quad&|\gamma \label{eq: UL}
			\end{align}
			\end{subequations}
		with $A=\nabla_y C(\tilde x,\tilde u,\tilde y)\in\mathbb{R}^{(8|\mathcal{N}|+4|\mathcal G|)\times|\mathcal{D}|}$, $b=C(\tilde x,\tilde u,\tilde y)\in\mathbb{R}^{(8|\mathcal{N}|+4|\mathcal G|)}$.
	 	As discussed in \cite{ralph2004some}, 
	 	penalty parameter $\rho$ can be modulated to meet the complementartiy satisfaction.



		\subsection{Local Convergence Analysis for Approximate BQP }
		In this subsection, we will show a convergence analysis of the simplified QP \eqref{eq:BQP}. For representational convenience, we introduce $\Delta y=y-\tilde{y}$ as the primal step.
		
		The merit function 
		\begin{equation}\label{eq: merit function}
		\psi(y)\doteq\frac{1}{2}y^\tr Q y +g^\tr y -\rho \phi(y)
		\end{equation}
		
		represents the outer loop objective function.
		It pointed out that merit
		function $\psi(y)$  at $y^k$ iteration is non-increasing towards Equation~\eqref{eq:BQP} for the local convexity of $\upsilon(y)$ \cite[Section III]{lcqp} and the property  
		\[\nabla\psi(y^k)^\tr \Delta y=\nabla\upsilon(y^k)^\tr \Delta y.\]
		However no convergence rate is discussed.  
		
		To ensure any local solution is a regular stationary point of \autoref{eq:BQP}, two assumptions are introduced below.
			\begin{assumption}
				\label{prop::LICQ}
					Linear Independence Constraint Qualification Condition (LICQ)\\
					The matrix $A$ has full row rank in the optimal value of $y^*$ thus the gradients of active inequality and equality constraints are linearly independent. We refer to \cite[Chapter 2]{Hall2021a} for further discussion.
			\end{assumption}
			\begin{assumption}
				\label{prop::SOSC}
				Second Order Sufficient Condition (SOSC)\\
				We assume the Hessian matrix $Q$ is negative semi-definite in a local neighborhood of $y^*$. This statement is well-known in Newton-type algorithms.  Suppose $Q$ is not negative semi-definite in the current iteration, then set $Q\leftarrow Q-\rho I$   (\cite{bliek1u2014solving}).
			\end{assumption}

			\begin{theorem}
				Let Assumption~\autoref{prop::LICQ} and \autoref{prop::SOSC} of Equation~\eqref{eq:BQP} be applicable, then the iteration of BQP can converge to the local saddle point $\phi(y^*)$ with super-linear convergence rate by using suitable line search step size $\alpha$ and penalty parameter $\rho$. 
			\end{theorem}
			\Proof
			The Lagrangian function of Problem~\eqref{eq:BQP} shows as
			\begin{equation}\label{eq:larg of BQP}
				\begin{split}
					\mathcal L(y,\lambda, \mu, \gamma)\doteq&\frac{1}{2}\Delta y^\tr Q \Delta y +(g-\rho \nabla\phi(\tilde y))^\tr \Delta y\\
					&+\lambda^\tr(b+A\Delta y
					)+\mu^\tr  y+\gamma^\tr(\mathbf{1}-y).
				\end{split}
			\end{equation}
			Assume $[\lambda_{\text{act}}^\tr, \mu_{\text{act}}^\tr, \gamma_{\text{act}}^\tr]^\tr$ collects the dual variables of the active inequalities of Problem~\eqref{eq:BQP}, $[(b+A\Delta y
			)_{\text{act}}^\tr, y_{\text{act}}^\tr, (\mathbf{1}-y)_{\text{act}}^\tr]^\tr$ collects the active constraints.
			The KKT (Karush-Kuhn-Tucker) system of Equation~\eqref{eq:larg of BQP} can be summarized as
			\begin{equation*}
				\begin{split}
					\begin{bmatrix}
						Q& A_{\text{act}}^\tr &I^{\mu}_{\text{act}} &-I^{\gamma}_{\text{act}}\\[0.16cm]
						A_{\text{act}}& 0& 0&0\\[0.16cm]
						I^{\mu}_{\text{act}}& 0&0 &0\\[0.16cm]
						-I^{\gamma}_{\text{act}}& 0&0 &0
					\end{bmatrix}\begin{bmatrix}
						\Delta y\\[0.16cm]
						\lambda_{\text{act}}\\[0.16cm]
						\mu_{_{\text{act}}}\\[0.16cm]
						\gamma_{\text{act}}
					\end{bmatrix}=\begin{bmatrix}
						-g+\rho \nabla\phi(\tilde y)\\[0.16cm]
						-b\\[0.16cm]
						0\\[0.16cm]
						0
					\end{bmatrix}.
				\end{split}
			\end{equation*}
			Thus the optimal updating primal variable can be expressed as  
			\[\Delta y^*(\rho)=-Q^{-1}\left(g-\rho \nabla\phi(\tilde y)+(A^\tr \lambda)_{\text{act}}+\mu_{\text{act}}-\gamma_{\text{act}}\right).\]
			
			Notice the residual norm of $\phi(y)$ can be expressed as
			\begin{equation}\label{eq:convergence}
				\begin{split}
					\left\|\frac{y^\tr \nabla\phi(\tilde y) }{\vspace{1mm}{\tilde  y^\tr \nabla\phi(\tilde y)}  }\right\|
					=&\left\|\frac{(\tilde y+\alpha\Delta y^*(\rho))^\tr \nabla\phi(\tilde y) }{\vspace{1mm}{\tilde y^\tr \nabla\phi(\tilde y)}  }\right\|\\[0.16cm]
					=&\left\|1+\alpha\frac{\Delta y^*(\rho)^\tr \nabla\phi(\tilde y) }{\vspace{1mm}{\tilde y^\tr \nabla\phi(\tilde y)}  }\right\|
				\end{split}
			\end{equation}
			with $\alpha$ denotes the step size. As long as
			\[\alpha=-\frac{\vspace{1mm}{\tilde y^\tr \nabla\phi(\tilde y)}}{\Delta y^*(\rho)^\tr \nabla\phi(\tilde y)},\]
			Equation~\eqref{eq:convergence} can provide a local super-linear convergence (\cite{nocedal2006numerical}) of the penalty function $\phi(y)$. 
			
	\hfill$\blacksquare$
			
			Note that this is an extension result of \cite[Section III-C]{lcqp} by considering the related active inequality. We give the corresponding local convergence result since the penalty parameter $\rho$ can also be turned.
			
			\subsection{AO-SBQP}
			 Integrated the structure of \emph{relaxed BQP},  
				\autoref{alg:AO-SBQP} summarizes the full AO-SBQP steps for solving problem \eqref{eq:object2}. The main idea is to use the alternate optimization method to solve the continuous and Boolean variables respectively.

			\begin{algorithm}[h!]
				\small
				\caption{AO-SBQP Method}\label{alg:AO-SBQP}
				\textbf{Input:}  initial guess $\tilde{y}$, initial $\bar x, \bar u$, a termination tolerance $\epsilon>0$, an initial factor $\rho>0$ and update rate $\beta>1$.\\
				\textbf{Repeat:} 
				\begin{enumerate}
					\item \emph{Linear Optimal Power Flow (AO1):}\label{Power flow step} Solve an optimization problem consists of \eqref{eq:object} and \eqref{eq:all equaltity} with given $\tilde{y}$. Then
					output optimal power system decision variables $\tilde{x}, \tilde u$.\footnotemark

					\item \emph{Sequential BQP (AO2):}
					\begin{enumerate}
						\item \emph{Globally Search:} Solve QP
						consists of
\eqref{eq:BQP}
with zero penalty parameter. Output optimal switch variable $\hat y$.\label{step: one of SBQP}
						\item \emph{Update Penalty Function Approximate:}
						\begin{equation*}
						\begin{split}
						\phi(y)&\approx \phi(\hat y) + (y-\hat y)^\top \nabla \phi(\hat y)\\
						&= (\phi(\hat y) -\hat y^\top \nabla \phi(\hat y))+y^\top \nabla\phi(\hat y).
						\end{split}
						\end{equation*}
						
						\item \emph{Locally Search:} Maximize the penalty QP \eqref{eq:BQP}.\label{sec:Locally Search}
						\item  \emph{Line Search and Inner Termination Criterion:} 
						
						$\alpha = StepLength(\hat y, \tilde{y}, \rho);$
						
						$\hat{y} \approx \hat{y} +\alpha (\tilde{y}-\hat y)$.
						
						Check if $|\phi(\hat y)|\leq \epsilon$, if not, go to Step \eqref{step:update}; if yes, $\tilde{y}\leftarrow\hat y$ and go to Step \eqref{step:outerterm}.
						\item \emph{Penalty Parameter Update:}\label{step:update}
						
						$\rho=\beta
						\cdot
						\rho$ 
						and return Step \eqref{step: one of SBQP}
					\end{enumerate}
					\item \emph{Outer Termination Criterion:}\label{step:outerterm}
					Check if \begin{equation*}
					\left\| ((\tilde{x},\tilde{u} )|\tilde{y})-(\bar{x},\bar{u})\right\|\leq \epsilon,
					\end{equation*}
					if not, go back to Step \eqref{Power flow step} and set $(\bar x, \bar u)=((\tilde{x},\tilde{u} )|\tilde{y}) $, $\bar y=\tilde{y}$; if yes, output the result.
				\end{enumerate}
				
				\textbf{Output:} $(x^*, u^*, y^*)\leftarrow (\tilde x,\tilde u,\tilde y)$.
			\end{algorithm}
			\footnotetext[1]{This step can be solved by any NLP solver.}
			\subsubsection{AO1}
			 Derived from  \cite{frank2012optimal1}, this variant of Optimal Power Flow aims to maximize 
			 the linear electric power
			 guarantee objective
			 with system input $u$, state $x$, fixed $\tilde{y}$ and the power grid physical constraints \eqref{eq::physical constyraints}.
			\subsubsection{AO2}
			The Hessian $Q$ and gradient $g$ are evaluated jointly with $(\tilde x,\tilde u,\tilde y,\tilde{\lambda})$. In order to search for a global maximizer of \autoref{eq:BQP} without the complementarity constraints, parameter $\rho$ of \autoref{eq:approximate} is set as $0$ in Step~\eqref{step: one of SBQP}. 
			 The later steps aiming to asymptotically meet \autoref{cons:lincompls} by increasing $\rho$ in Step~\eqref{step:update}. Note that, both Step~\eqref{step: one of SBQP} and Step~\eqref{sec:Locally Search} are simple convex QP which can be solved by any stable QP solver. We refer \cite[Section III]{2022Alternating} as a reference for more details.
			\begin{remark} \emph{(Relaxation of AO2)}\label{remark 1}
			For the robustness of switching from AO2 to AO1, \eqref{eq:BQP:affcons} can be arbitrarily replaced by the following Inequality~\eqref{eq: relaxed constraints}
			in Step~\eqref{step: one of SBQP} and Step~\eqref{sec:Locally Search} which named as \emph{mixed AO2},
			\begin{equation}\label{eq: relaxed constraints}
			\left\{			\begin{split}
			&\sum_{k\in\mathcal{G}}p_k^g-\sum_{k\in\mathcal{D}}y_kp_k^d\geq 0\\
			&\sum_{k\in\mathcal{G}}\overline{q_k^g}-\sum_{k\in\mathcal{D}}y_kq_k^d\geq 0\\
			&\sum_{k\in\mathcal{D}}y_kq_k^d-\sum_{k\in\mathcal{G}}\underline{q_k^g}\geq 0.
			\end{split}\right.
			\end{equation}
			Or even, in some cases, the entire AO2 process can be replaced by solving the \emph{relaxed AO2} \eqref{eq: preprocess} module below,
			\begin{equation}\label{eq: preprocess}
			\begin{split} 
			\max_{y}\;\;\sum_{k\in\mathcal{D}}y_k^2 r_k  p^d_k  \quad\mathrm{s.t.}\;\; \eqref{eq: relaxed constraints}, \eqref{eq: MPCC}
			\end{split}.
			\end{equation}
			\end{remark}

			\section{Numerical Result}\label{sec: numerical result}
			
			In this section, we illustrate the numerical result of AO-SBQP method drawing upon the modified 30-bus power network. The power sources of microgrid are considered into the system.
			
			\subsection{Data and Implementation}
			The problem data is obtained from \texttt{MATPOWER} dataset \cite{Zimmerman2011} and the implementation of Algorithm~\eqref{alg:AO-SBQP} relays on \texttt{Casadi-v3.5.5} with \texttt{IPOPT} 
			(\cite{casadi}). Though \texttt{MATPOWER} repository is for transmission networks which are high voltage networks, the mathematical model of microgrid is same.
			
In the modiﬁed 30-bus case, $\mathcal G=\{1,2,13,22,23,27\}$. To create a \emph{demand-to-power mismatch scenario}, we increase the active and reactive power demands of \newshijie{buses with load} by $2.5\; \text{p.u.}$ (per unit) and $0.7\; \text{p.u.}$ respectively. The lower and upper bounds of reactive power inputs are reduced to half of the previous ones while \newshijie{upper bound of active ones are reduced to 70\%} . In addition, $r_k$'s are randomly set into five levels from 1 to 5 for each demand and the criterion (i.e. complementarity tolerance) is set as
$10^{-6}$. Problem~\eqref{eq:object2} consists $60$ status, $12$ system inputs and $30$ switch variables \newshijie{when all buses contain consumers}.
				\begin{figure}[H] 
				\centering
				\includegraphics[width=0.6\linewidth,height=0.25\textheight]{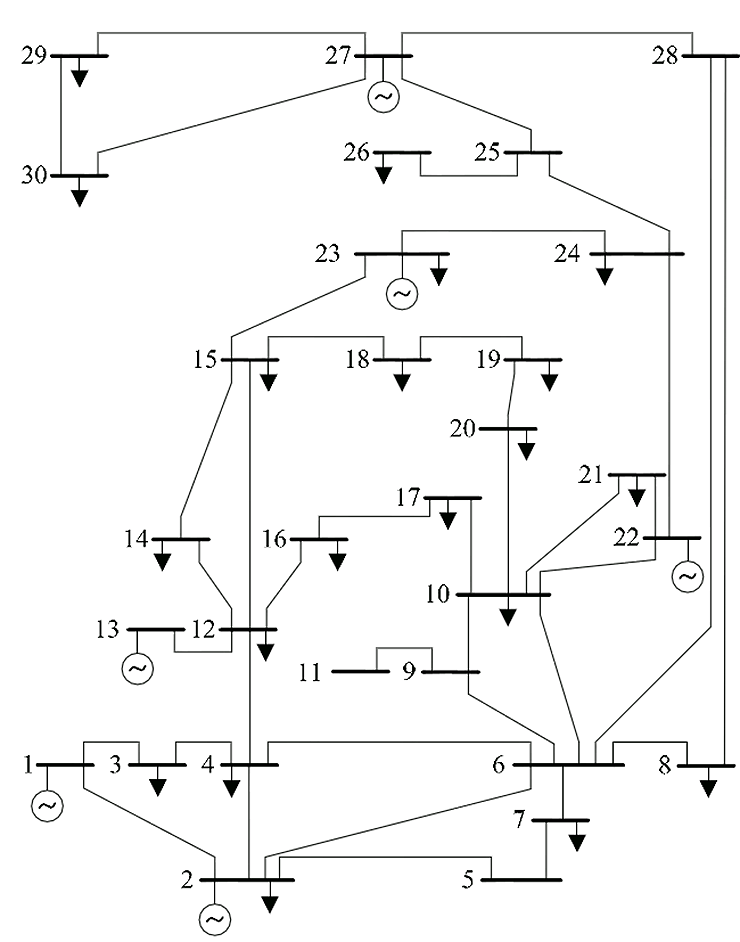}
				\caption{Modiﬁed IEEE 30-bus system from \cite{christie2000power}. } 
				\label{fig:ieee30bus}
			\end{figure}
			
Since we did not set multi demands for each bus, this indicates $|\mathcal D|=|\mathcal N|$, and \eqref{eq:BQP:affcons} is going to be an overdetermined system. Therefore our implementation focus on the other three relax versions of AO2 mentioned in Remark~\ref{remark 1}.

			\subsection{Numerical Comparison}
\begin{figure}[H] 
				\centering
				\includegraphics[width=0.7\linewidth,height=0.2\textheight]{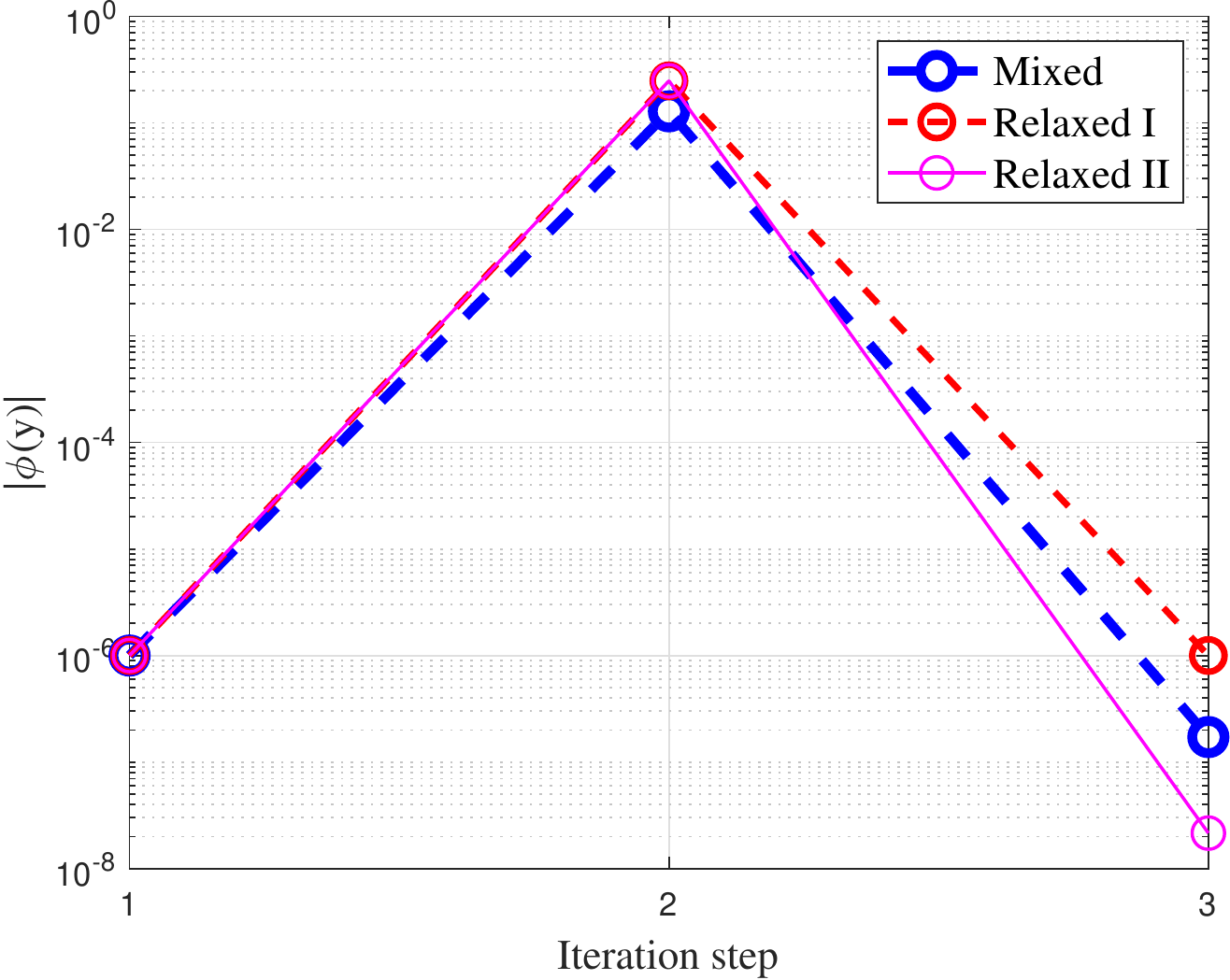}
				\caption{Convergence of $\phi(y)$ with three
				 variations of AO2.} 
				\label{fig:phiconvergence}
			\end{figure}

	In this section, we show the comparison of Algorithm~\eqref{alg:AO-SBQP} with different variations of AO2. Note that all the variations consist inequality constraints $\eqref{eq: relaxed constraints}, \eqref{eq: BL}, \eqref{eq: UL}$, and the only difference are the objectives, 
	a) Mixed:  \eqref{eq:BQP:obj}, b) Relaxed I: $\sum_{k\in\mathcal{D}}y_k^2 r_k  p^d_k-\rho \phi(y)$ , c) Relaxed II: $\sum_{k\in\mathcal{D}}y_k^2 r_k  p^d_k -\rho \nabla\phi(\tilde y)^\tr y$.
%
%
%
				\begin{figure}[H] 
				\centering
				\includegraphics[width=0.7\linewidth,height=0.2\textheight]{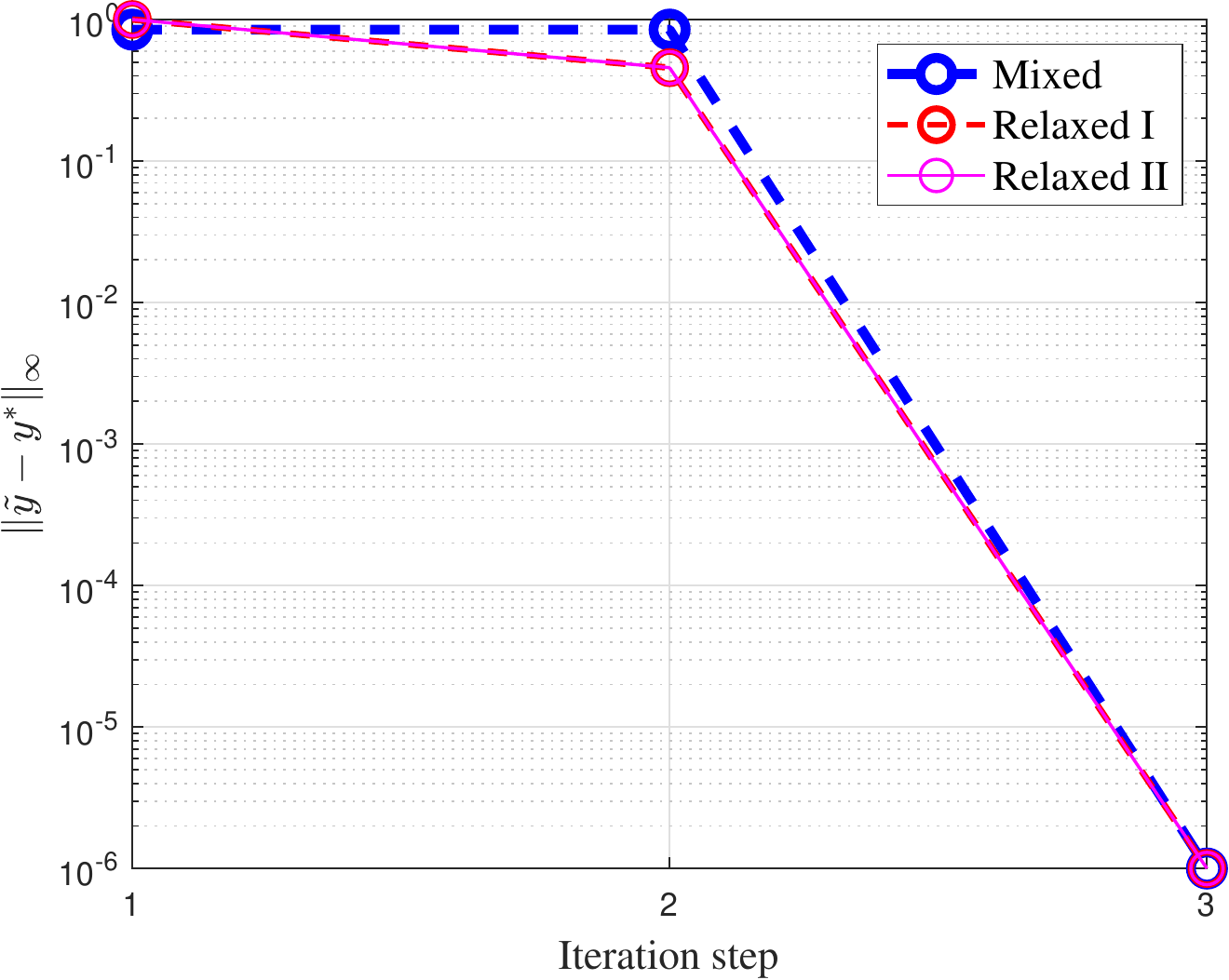}
				\caption{Convergence of $y$ with three variations of AO2. } 
				\label{fig:yconvergence}
			\end{figure}	
				Figure~\eqref{fig:yconvergence} and Figure~\eqref{fig:phiconvergence} shows the convergence of the switch variable $y$ and the complementarity satisfaction $\phi(y)$ respectively. It can be seen that even if the dimension of $y$ is $30$, all the variations of SBQP can converge to the given complementarity tolerance in only a few steps but converge to different solutions. This shows completely different properties than the B$\&$B based solvers.	
			\begin{table}[h]
				\caption{Convergence Comparison}
				\label{table:comparison1}
				\begin{center}
					\begin{tabular}{|c||c|c|c|}
						\hline
					Method	&\makecell{Mixed}& \makecell{Relaxed I}&\makecell{Relaxed II}\\
						\hline
						$|\phi(y)|$ &\makecell{1.7e-07}  &\makecell{1.0e-06} &\makecell{2.2e-08}\\
						\hline
						Time(s)  &\makecell{0.021} &\makecell{0.046} &\makecell{0.019}\\
						\hline
						Iteration &\makecell{2}  & \makecell{2}&\makecell{2}\\
						\hline
					\end{tabular}
				\end{center}
			\end{table}
		Table~\ref{table:comparison1} shows the  comparison among complementarity satisfaction, operation time and iteration by using different variations. 
		It can be seen that all indicators of the three methods are similar which is different from the results seen in \cite{2022Alternating}. Since the  
	    three relaxation variations in this paper decrease the number of inequality constraints in AO2 that induce it easier to solve.
	    Notice that, due to the nonlinear structure of \eqref{eq:object}, the implementation is more complex than \eqref{eq: original obj}, therefore only AO1 benefits from the reformulation of \eqref{eq:object}. 
	    
			\begin{table}[h]
				\renewcommand{\arraystretch}{1.2}
				\caption{Performance Comparison}
				\label{table:comparison2}
				\begin{center}
						\begin{tabular}{|c||c|c|c|}
						\hline
					Method	&\makecell{Mixed}& \makecell{Relaxed I}&\makecell{Relaxed II}\\
						\hline
						Objective   &\makecell{5.421}  &\makecell{2.428} &5.877\\
						\hline
						$\sum\limits_{k\in \mathcal D}y_kp_k^d$ (p.u.)    &\makecell{1.567}  &\makecell{0.825} &\makecell{1.592}\\
						\hline
						$\sum\limits_{k\in \mathcal D}y_kq_k^d$ (p.u.)   &0.702  &0.347 &0.709\\
						\hline
					\end{tabular}
				\end{center}
			\end{table}
		Table~\ref{table:comparison2}
		shows the comparison of final performance comparison by using different methods. As can be seen, different variations converge to different local optimum. At least in this case, Relaxed II gets a bit better performance than the other two and Relaxed I is a bit conservative. Note that, both Table~\ref{table:comparison1} and Table~\ref{table:comparison2} show the benefits of the approximate BQP method (Mixed and Relaxed II).

		\begin{table}[h]
			\renewcommand{\arraystretch}{1.2}
			\caption{Optimal Power Inputs}
			\label{table:optimal inputs}
			\begin{center}
			\begin{tabular}{|c||c|c|c|c|c|c|}
					\hline
					Bus $\#$& 1&2&13&22&23&27\\
					\hline
				pg (p.u.)
					&0.400
					&0.400
					&0.134
					&0.250
					&0.150
					&0.275\\
					\hline
				qg (p.u.) 
					&-0.037
					&0.161
					&0.213
					&0.246
					&0.070
					&0.125
					 \\
					\hline
				\end{tabular}
			\end{center}
		\end{table}
			Table~\eqref{table:optimal inputs} shows the numerical result of optimal system inputs by using
			Relaxed II based Algorithm~\eqref{alg:AO-SBQP}. The active demands relate with bus $\{$1,2,4,9,10,11,12,13,14,16,17,18,\\21,22,23,24,25,26,30$\}$.

			\section{Conclusion}
			
			In this paper, we proposed
			an effective and fast convergence method named AO-SBQP to optimize microgrid demand shut-offs problems. Importantly, local convergence theory of approximate BQP has been proposed.
			Moreover, a numerical result on modiﬁed IEEE 30-bus case study illustrates the potential of AO-SBQP in this area. 
			Different from B$\&$B and SDR,
			AO-SBQP can achieve a feasible local optimal solution without 
			tree storage structure
			or
			matrix variables.   
			Future research will investigate multistage optimal demand shut-offs and time varing priority of single bus-multiple demands on larger case studies. Moreover, the rank evaluating priority of each demand can vary with time.
			Comparison of accuracy and computation time of B$\&$B and SDR will also be considered.

			\bibliography{paper,zsjcitation} 
		\end{document}